\pdfoutput=1

\RequirePackage{fix-cm}

\documentclass[leqno, fleqn, centertags, 12pt]{article}

\usepackage{latexsym}
\usepackage{amsmath}
\usepackage{amssymb}
\usepackage{verbatim}

\usepackage[T1]{fontenc}

\usepackage[upright, widespace]{fourier}

\usepackage{bm}

\usepackage{fullpage}

\makeatletter
\renewcommand{\@makefntext}[1]{\vspace*{0.5ex}\parindent=0em
\hspace*{-0.4em}
\hbox to 0.4em{\hss\@makefnmark}\hspace*{0.4em}{#1}
}
\makeatother

\newcommand{\myuppar}[1]{\vspace{\medskipamount}\textbf{#1}\hspace*{0.5em}}

\newcounter{myparnum}
\renewcommand{\themyparnum}{
\arabic{myparnum}}
\newcommand{\mypar}[2]{\refstepcounter{myparnum}{\vspace{\medskipamount}\textbf{{\themyparnum. #1}\label{#2}}\hspace{0.5em}}}

\newcommand{\proof}{\vspace{\medskipamount}{\textbf{{\emph{Proof}.}}\hspace*{1em}}}
\newcommand{\eproof}{ $\blacksquare$}

\def\sss{\hspace{0.05em}\ }
\def\dss{\hspace{0.1em}\ }

\def\qss{\hspace{0.2em}\ }

\def\oss{\hspace{0.4em}\ }

\def\halfff{\hspace*{0.025em}}
\def\fff{\hspace*{0.05em}}
\def\dff{\hspace*{0.1em}}
\def\trf{\hspace*{0.15em}}
\def\qff{\hspace*{0.2em}}

\newcommand{\nsp}{\hspace*{-0.1em}}
\newcommand{\nnsp}{\hspace*{-0.15em}}

\newcommand{\dnsp}{\hspace*{-0.2em}}

\renewcommand{\leq}{\leqslant}
\renewcommand{\geq}{\geqslant}

\newcommand{\nnn}{\mathbf{N}}

\newcommand{\id}{\operatorname{id}}
\newcommand{\norm}[1]{\|\qff #1 \qff\|}
\newcommand{\sco}[1]{\langle\trf #1 \trf\rangle}

\newcommand{\ttoo}{\hspace*{0.2em}\longrightarrow\hspace*{0.2em}}

\begin{document}

\setlength{\baselineskip}{12pt plus 0pt minus 0pt}
\setlength{\parskip}{12pt plus 0pt minus 0pt}
\setlength{\abovedisplayskip}{12pt plus 0pt minus 0pt}
\setlength{\belowdisplayskip}{12pt plus 0pt minus 0pt}

\newskip\smallskipamount \smallskipamount=3pt plus 0pt minus 0pt
\newskip\medskipamount   \medskipamount  =6pt plus 0pt minus 0pt
\newskip\bigskipamount   \bigskipamount =12pt plus 0pt minus 0pt

\author{Nikolai\qss V.\qss Ivanov\quad and\quad Marina\qss Prokhorova}
\title{The unitary group in the strong topology\\
and a construction of\dss Dixmier--Douady} 
\date{}

\footnotetext{\hspace*{-0.65em}\copyright\qss 
Nikolai\qss V.\qss Ivanov\qss and\qss Marina\qss Prokhorova,\qss 2025.\oss 
}

\footnotetext{\hspace*{-0.65em}The work of M. Prokhorova 
was done during a postdoctoral fellowship at the University of Haifa
and was partially supported by the ISF grant 876/20.\qss }

\maketitle

\vspace*{3ex}

\renewcommand{\baselinestretch}{1.14}
\selectfont

\myuppar{Introduction.}
For a separable infinite-dimensional\dss Hilbert space $H$ let $U(H)$ be its unitary
group equipped with the norm topology and\dss let $\mathcal{U}(H)$ be the same
group equipped with the strong topology.\qss
By a well known theorem of\dss Kuiper\dss the space $U(H)$ is contractible.\qss
By a theorem of\dss Dixmier--Douady\dss \cite{dd} 
the space $\mathcal{U}(H)$ is also contractible.\qss
See\qss \cite{dd},\dss Lemma\qss 3.\qss
The proof\sss of\dss Dixmier--Douady\dss is short and striking.\qss 
It is based on the realization of $H$ as the space $L_2([0,\dff 1])$
and an explicit construction of families of subspaces and operators in $H$ 
with rather special\sss properties.\qss 
See\qss \cite{dd},\qss Lemma\qss 2,\qss
or\sss the section\dss \emph{Dixmier--Douady--like families}\qss below.\qss
Unfortunately,\qss this proof\sss leaves hidden\sss 
the geometric meaning of\dss the theorem.\qss
The first goal of this note is\sss to give a direct geometric proof of\dss
the contractibility of $\mathcal{U}(H)$\nnsp.

Dixmier--Douady construction can also be used to prove the contractibility
of other spaces of operators.\qss
See,\qss for example,\qss Atiyah--Segal\qss \cite{as},\qss Proposition\qss A2.1.\qss
Recently\sss it was used by the second author\qss \cite{p}.\qss
See\qss \cite{p},\qss Propositions\qss 4.1\qss and\qss 4.3.\qss
This motivated\sss the authors to
provide a geometric analogue
of\dss Dixmier--Douady\dss construction.\qss
Namely,\qss we will\sss give a geometric construction,\qss 
and this is the second goal of\sss this note,\qss 
of families of subspaces and operators which 
we call\sss  \emph{Dixmier--Douady--like\dss families}.\qss
They have all\sss properties of\dss Dixmier--Douady\dss families
essential\sss for applications,\qss
and have better continuity properties.\qss
This geometric construction sheds additional\dss light on\dss Dixmier--Douady\dss results
and\sss their applications.\qss

Let\sss us review\sss the basic properties of\sss $\mathcal{U}(H)$\nnsp,\qss
although\sss they play only a marginal\sss role in our arguments.\qss
They are conveniently discussed in\dss M.\dss Schottenloher's\dss paper\dss \cite{s}.\qss
First of all,\dss $\mathcal{U}(H)$ is a topological\sss group.\qss
See\qss \cite{s},\qss Proposition 1.\dss
The space $\mathcal{U}(H)$ is metrizable and\sss hence compactly generated\dss
(we will not use this fact at all).\qss
See\qss \cite{s},\qss Proposition 3 (note that $H$ is assumed to be separable).\qss
The topology of\sss $\mathcal{U}(H)$ (i.e.\dss the strong topology)
coincides with the compact-open topology.\qss
See\qss \cite{s},\qss Proposition 2.\qss
The same argument\sss proves a more general\sss fact,\qss
which is useful\sss for applications.\qss
Namely,\qss for every $C>0$ the strong topology
coincides with the compact-open topology on 
the space of operators in $H$ with the norm $\leq C$\nnsp.\qss

The following\sss lemma\sss is essentially all\sss what we need to
know about the convergence in the strong topology,\qss
beyond the fact that it is implied by the convergence in the norm topology.\qss

\mypar{Lemma.}{l10}
\emph{Let\dss
$H_1\subset\dff H_2\subset\dff H_3\subset \ldots$
be a sequence of subspaces of\dss $H$ such that\sss
their union is dense in $H$\nnsp.\qss
If\dss $f_t\fff,\dff t\in [0,\dff 1)$
is a family of\sss operators in $H$ with the norm $\leq 1$
such that $f_t$ is equal\sss to the identity on $H_i$
for $t\geq 1-1/i$\nnsp,
then\sss $f_t$ strongly converges to $\id_{\dff H}$ when $t\ttoo 1$\nnsp.\qss
More generally, let
$X$ be a topological space
and\dss $f_{x\fff,\dff t}\dff,\qff x\in X\fff,\qff t\in [0\fff,\dff 1)$
be a family of\sss operators in $H$ with the norm $\leq 1$
such that $f_{x\fff,\dff t}$ is equal\sss to the identity on $H_i$
for $t\geq 1-1/i$\nnsp.\qss
Let $f_{x\fff,\dff 1}=\id_{\dff H}$ for $x\in X$\nnsp.\qss
Then the family $f_{x\fff,\dff t}\dff,\qff x\in X\fff,\qff t\in [0\fff,\dff 1]$
is strongly continuous at $X\times 1$\nnsp.}

\proof
Let us prove the first claim.\qss
It\sss is sufficient to prove that $f_t(h)$ converges to $h$
in the norm when $t\ttoo1$ 
for every $h\in H$\nnsp.\qss
By the assumption,\qss for every $\varepsilon>0$ there exist $i$ and $h'\in H_i$
such that $\norm{h-h'}<\varepsilon/2$\nnsp.\qss
Then $f_t(h')=h'$ for $t\geq 1-1/i$ and for such $t$ we have
\[
\quad
\norm{f_t(h)-h}
\qff \leq\qff
\norm{f_t(h)-f_t(h')}\qff +\qff \norm{f_t(h')-h'}\qff +\qff \norm{h'-h}
\]

\vspace{-39pt}
\[
\quad
\phantom{\norm{f_t(h)-h}
\qff }=\qff
\norm{f_t(h-h')}\qff +\qff \norm{h'-h'}\qff +\qff \norm{h'-h}
\]

\vspace{-39pt}
\[
\quad
\phantom{\norm{f_t(h)-h}
\qff }\leq\qff
\norm{h-h'}\qff +\qff 0\qff +\qff \norm{h'-h}
\qff =\qff
2\dff\norm{h'-h}
\qff <\qff
\varepsilon
\qff.
\]
This proves that $f_t(h)$ converges to $h$ and hence proves the first statement.\oss  

The estimates in the previous paragraph depend on the properties of $f_t$\nsp,\qss
but not on $f_t$\nsp.\qss
Hence we can replace there $f_t$ by $f_{x\fff,\dff t}$ and conclude
that with the same choice of $i$ we have
$\norm{f_{x\fff,\dff t}(h)-h}
\qff <\qff
\varepsilon$
for $t\geq 1-1/i$ and every $x\in X$\nnsp.\qss
The second claim follows.\oss \eproof\vspace{-1pt}

\mypar{Theorem\dss (Kakutani).}{t11}
\emph{The unit sphere in $H$ is contractible.\oss}\vspace{-1pt}

\proof
See Kakutani\qss \cite{k}\qss for a beautiful elementary proof.\oss  \eproof\vspace{-1pt}

\myuppar{Infinitely dimensional\dss Stiefel\dss manifolds.}
Let \nsp$V_k(H)$\nsp\ be the space of $k${\dnsp}-tuples 
$(v_1\fff,\dff v_2\fff,\dff \ldots\fff,\dff v_k)$\nnsp,\qss
where $v_i\in H$ are of\sss length $1$ and pair-wise orthogonal.\oss
The topology of $V_k(H)$ is induced from the product of $k$ copies of $H$\nnsp.\qss
For $x\in V_k(H)$ let $H_x$ be the orthogonal complement of\sss the span
of\sss the vectors from the $k${\dnsp}-tuple $x$\nnsp.\qss
Let $p_x$ be the orthogonal\sss projection $H\to H_x$\nsp,\qss
and $i_x\colon H_x\to H$ be the inclusion.\qss
Let $G_k(H)$ be the space of pairs $(x,\dff h)$ such that $x\in V_k(H)$
and $h\in H_x$\nsp.\qss
Let $\pi_k\colon G_k(H)\ttoo V_k(H)$ be the map
$(x,\dff h)\longmapsto x$\nnsp.\qss
For most of the readers the following lemma will\sss be either obvious or known.\qss

\mypar{Lemma.}{l12}
\emph{The map\sss $\pi_k$
is a locally trivial\sss bundle with the structure group $U(H)$\nnsp.}

\proof
Clearly,\dss the maps $p_x$ and $i_x$ are adjoint to each other,
$i_x^*=p_x$\nsp.\qss
Since the target of the map $p_x$ depends on $x$\nnsp,\dss
we cannot directly speak about the continuity of the map $x\longmapsto p_x$\nsp.\qss
But the composition $i_x\dff p_x$ is a map $H\ttoo H$
and continuously depends on $x$ in the norm topology,\qss
as one can easily check.\qss

Let us fix a point $a\in V_k(H)$\nnsp.\qss
For $x\in V_k(H)$ belonging to a neighborhood $N$ of $a$ the maps
\[
\quad
r_x = p_a \dff |\dff H_x
=p_a\dff i_x
\colon H_x\ttoo H_a
\quad
\mbox{and}\quad
s_x = p_x \dff |\dff H_a
=p_x\dff i_a
\colon H_a\ttoo H_x
\qff
\]
are linear isomorphisms.\qss
The polar decomposition shows that 
$u_x
\qff =\qff
s_x\dff (s_x^*\dff s_x)^{\dff -1/2}$
is a unitary (i.e. isometric) isomorphism $H_a\ttoo H_x$\nsp.\qss
Since the operator $s_x\dff (s_x^*\dff s_x)^{\dff -1/2}$ is unitary,\qss
its inverse $u_x^{\dff -1}$ is equal\sss to its adjoint,\dss i.e.\qss to 
$u_x^*
\qff =\qff 
(s_x^*\dff s_x)^{\dff -1/2}\dff s_x^*$\nsp.\qss

At the same time 
if\sss $w\in H_a$ and $v\in H_x$\nsp,\oss then
$\sco{p_a(v),\dff w}=\sco{v,\dff w}=\sco{v,\dff p_x(w)}$
and hence the operators $r_x$ and $s_x$ are adjoint to each other,\dss
$s_x^*=r_x$\nsp.\qss
It follows that
\[
\quad
s_x^*\dff s_x
\qff =\qff
r_x\dff s_x
\qff =\qff
p_a\dff i_x\dff p_x\dff |\dff H_a
\qff.
\]
Since $i_x\dff p_x$ continuously depends on $x$ in the norm topology,\qss
this implies that $s_x^*\dff s_x$ continuously depends on $x$ (in the norm topology).\qss
In turn, this implies that $(s_x^*\dff s_x)^{\dff -1/2}$ and 
$i_x\dff u_x
\qff =\qff 
i_x\dff s_x\dff (s_x^*\dff s_x)^{\dff -1/2}
\qff =\qff 
i_x\dff p_x\dff (s_x^*\dff s_x)^{\dff -1/2}$ 
continuously depend on $x$\nnsp.\qss
As above,\qss we cannot directly speak about the continuity of 
$x\longmapsto u_x$\nsp.\qss
Similarly,\qss in order to discuss the continuity of  
$x\longmapsto u_x^{\dff -1}\qff =\qff u_x^*$ on $x$\nnsp,\qss
we need\sss to compose $u_x^*$ with $p_x$\nnsp.\qss
This composition is equal\sss to
\[
\quad
u_x^*\dff p_x
\qff =\qff
u_x^*\dff i_x^*
\qff =\qff
(i_x\dff u_x)^*
\qff
\]
and hence continuously depends on $x$\nnsp.\qss

The map $N\times H_a\ttoo  \pi_k^{\dff -1}(N)$
given by $(a,\dff h)\longmapsto (x,\dff u_x(h))$ is,\qss clearly,\qss a bijection.\qss
Since $i_x\dff u_x$ continuously depends on $x$ in the norm topology
and the topology of $G_k(H)$ is induced from $V_k(H)\times H$\nnsp,\qss
this map is continuous.\qss
The inverse map is given by
\[
\quad
(x,\dff h)\longmapsto (a,\dff u_x^{\dff -1}(h))
\qff =\qff
(a,\dff u_x^*\dff p_x(h))
\qff.
\]
Since $u_x^*\dff p_x$ continuously depends on $x$\nnsp,\qss
this map is also continuous.\qss
It\sss follows that the above map and its inverse provide
a local\sss trivialization of $\pi_k$ near $a$\nnsp.\qss
Moreover,\qss such trivializations for different $a$ and $N$
differ by norm continuous maps.\qss
The lemma follows.\oss  \eproof

\myuppar{Remarks.}
The above proof\sss is an adaptation of some arguments of\qss Dixmier--Douady\qss \cite{dd}.\qss
See\qss \cite{dd},\qss the proof of\dss Theorem\dss 2\dss and the subsection\dss 9.\qss
We will\sss need only the case $k=1$ of this lemma,\qss
but the proof in this case is no simpler than in general.\qss

\mypar{Theorem.}{t13}
\emph{\dnsp$\mathcal{U}(H)$ is contractible.\qss
Moreover,\qss
the contracting homotopy can be chosen to be continuous also as a homotopy of maps 
$U(H)\ttoo U(H)$\nnsp,\qss
except the last moment.\qss}

\proof
Let us consider the bundle
$\pi_1\colon G_1(H)\ttoo V_1(H)$ and identity
the base $V_1(H)$ of\sss this bundle with the unit sphere in $H$\nnsp.\qss
The unit sphere is metrizable and\sss hence paracompact.\qss
By\dss Lemma\qss \ref{l12}\qss the bundle $\pi_1$ is\sss locally\sss trivial.\qss
Also,\qss the base is contractible by\dss Theorem\dss \ref{t11}.\qss
It\sss follows that this bundle is\sss trivial.\qss
Let us fix a trivialization of\sss this bundle.\qss
This trivialization defines for every $x,\dff y\in V_1(H)$ a unitary isomorphism
$h(y, x)\colon H_x\ttoo H_y$ such that the following properties hold:\dss
$h(x, x)\qff =\qff \id$ for every $x\in V_1(H)$\nnsp;\dss
$h(z, y)\circ h(y, x)\qff =\qff h(z, x)$ for every  $x,\dff y,\dff z\in V_1(H)$\nnsp;\dss
and $i_y\dff h(y, x)\dff p_x$ continuously depends on $x,\dff y$ in the norm topology.\qss

Let us fix an orthonormal\sss basis $e_1\fff,\dff e_2\fff,\dff \ldots$ of\dss $H$\nnsp.\qss
For every unitary operator $u\colon H\ttoo H$ we will construct a family of unitary operators
$u_t\colon H\ttoo H$\nnsp, $t\in [0,\dff 1)$ such that $u_{\dff 0}=u$ and
$u_t(e_i)=e_i$ for every $u$ and $t\geq i/(i+1)=1-1/(i+1)$\nnsp.\qss
Moreover,\qss the map
$(u, t)\longmapsto u_t$
will\sss be continuous as map
$\mathcal{U}(H)\times [0,\dff 1)
\ttoo
\mathcal{U}(H)$
and also as a map
$U(H)\times [0,\dff 1)
\ttoo
U(H)$\nnsp.\qss
We will construct such operators in stages numbered by  $i\in\nnn$\nnsp.\qss

At the first stage we consider the action of
unitary 
operators $u\colon H\ttoo H$ on the vector $e_1$\nsp.\qss
Let us fix a contracting homotopy of $V_1(H)$ to $e_1$ 
and write it in the form
$x\longmapsto x_t$\nnsp,\dss $t\in [0,\dff 1/2]$\nnsp,\qss
so that $x_{\dff 0}\qff =\qff x$ and $x_{1/2}\qff =\qff e_1$ for every $x\in V_1(H)$\nnsp.\qss
For an unitary operator $u\colon H\ttoo H$ and $t\in [0,\dff 1/2]$ let $u_t\colon H\longrightarrow H$
be the unitary operator taking $e_1$ to $u(e_1)_t$ and acting as
$h(u(e_1)_t\dff,\dff u(e_1))\circ u$ on $H_{\dff e_1}$\nsp.\qss
Then $u_{\dff 0}\qff =\qff u$ and $u_{1/2}(e_1)\qff =\qff e_1$\nsp.\qss
The rule $(u, t)\longmapsto u_t$ defines maps
$\mathcal{U}(H)\times [0,\dff 1/2]
\ttoo
\mathcal{U}(H)$
and
$U(H)\times [0,\dff 1/2]
\ttoo
U(H)$\nnsp.\qss
Both are continuous because 
$u(e_1)$ continuously depends on $u$ in
strong and norm topologies and 
$i_y\dff h(y, x)\dff p_x$ continuously depends on $x,\dff y$ in the norm topology.\qss

At the second stage we apply
the same argument to the orthogonal complement $H_{\dff e_1}$ of $e_1$
and the vector $e_2\in H_{\dff e_1}$\nsp.\qss
For each unitary operator $u\colon H\ttoo H$ and $t\in [1/2,\dff 3/4]$ 
we get a unitary operator $u_t\colon H\ttoo H$
such that $u_t(e_1) = e_1$ for every $t\in [1/2,\dff 3/4]$ and 
$u_{\dff 3/4}(e_2)=e_2$\nsp. 
By continuing in this way we get the promised operators $u_t$\nsp, $t\in [0,\dff 1)$\nnsp.\qss
Let us set $u_1= \id_{\dff H}$ for every $u$\nnsp.\qss
Lemma\qss \ref{l10}\qss applied to $X=\mathcal{U}(H)$ 
and $f_{u\fff,\dff t}=u_t$ implies that the map
$(u, t)\longmapsto u_t$ is continuous
as a map
$\mathcal{U}(H)\times [0,\dff 1]\ttoo \mathcal{U}(H)$\nnsp.\qss
This proves the first statement of the theorem.\qss
The previous paragraph implies that the map
$U(H)\times [0,\dff 1)
\longrightarrow
U(H)$ 
given by the same rule is also continuous,\qss proving the second statement.\oss  \eproof

\myuppar{Remarks.}
The restriction to  $\mathcal{U}(H)\times [0,\dff 1)$ 
of the contracting homotopy constructed in the above proof\sss has the form
$(u,\dff t)\longmapsto h(u,\dff t)\circ u$ for a map
$h\colon \mathcal{U}(H)\times [0,\dff 1)\ttoo \mathcal{U}(H)$\nnsp.\qss
Moreover,\qss by the construction $h$ is continuous as a map
$\mathcal{U}(H)\times [0,\dff 1)\ttoo U(H)$\nnsp.\qss

The proof\sss of\dss Theorem\qss \ref{t13}\qss 
uses only that the operators $u$\sss are isometric embeddings\qss
(this is needed at the second stage and\sss later).\qss
Hence it
shows that the space of\sss isometric embeddings
$H\ttoo H$ with the strong topology\sss is\sss contractible
by a homotopy such as above.\qss

\myuppar{Dixmier--Douady--like families.}
We will construct a family $H_t\fff,\qff t\in [0,\dff 1]$
of subspaces of $H$ such that $H_{\fff 0}=\dff 0$ and $H_{\fff 1}=\dff H$\nnsp.\qss
Let $P_t\colon H\ttoo H_t$ be the orthogonal projection onto $H_t$ and let
$I_t\colon H_t\ttoo H$ be the inclusion.\qss
The operator $I_t\dff P_t\colon H\ttoo H$ will continuously depend on $t$ in the strong topology.\qss
In contrast with\dss Dixmier--Douady,  $I_t\dff P_t$ will
continuously depend on $t$ in the norm topology for $t\in (0,\dff 1)$\nnsp.\qss
Of course,\qss $I_t\dff P_t$ 
cannot continuously depend on $t$ in the norm topology for all $t\in [0,\dff 1]$\nnsp.\qss
The continuity in the strong topology will\sss follow from Lemma\dss \ref{l10}.\qss
Let $H'_t\qff =\qff H\ominus H_t$ be the orthogonal complement of $H_t$\nsp,\dss
let $Q_t\colon H\ttoo H'_t$ be the orthogonal projection,\qss and
$I'_t\colon H'_t\ttoo H$ be the inclusion.\qss
Then $I'_t\dff Q_t\qff =\qff \id_{\dff H}\dff -\qff I_t\dff P_t$\nsp.\qss

In addition,\qss we will construct\sss families of\sss 
isometric isomorphisms $U_t\colon H_t\ttoo H$\nnsp,\dss
$t\in(0,\dff 1]$
and $V_t\colon H'_t\ttoo H$\nnsp,\dss $t\in [0,\dff 1)$
such that $U_1=V_0=\id_{\dff H}$ and\sss the operators
$U_t\dff P_t$ and $V_t\dff Q_t$ 
norm continuously depend on $t$ for $t\in (0,\dff 1)$
and strongly continuously depend on $t$ at $t=0$ and $t=1$ respectively.\qss
Moreover,\qss the operators $I_t\dff U_t^{\dff -1}$ and $I'_t\trf V_t^{\dff -1}$
will strongly continuously depend on $t$ at $t=1$ and $t=0$ respectively.\qss
Note that\sss 
$I_{\dff 1}\dff U_1^{\dff -1}
\qff =\qff
I'_0\trf V_0^{\dff -1}
\qff =\qff
\id_{\dff H}$\nsp.\qss 

For such families $U_t\dff P_t$ strongly converges to $0$ when $t\ttoo 0$\nnsp,\qss
and  $V_t\dff Q_t$ strongly converges to $0$ when $t\ttoo 1$\nnsp.\qss
Indeed, $\norm{U_t\dff P_t(f)}=\norm{P_t(f)}\ttoo 0$ when $t\ttoo 0$ for
every $f$ because $I_t\dff P_t$ strongly converges to $0$\nnsp.\dss
The claim about $V_t\dff Q_t$ is proved in the same way.\qss
Also,\dss  $I_t\dff U_t^{\dff -1}$ and $I'_t\dff V_t^{\dff -1}$
norm continuously depend on $t$
for $t\in (0,\dff 1)$\nnsp.\qss
Indeed,\dss $I_t\dff =\dff P_t^*$ 
and\sss hence
$I_t\dff U_t^{\dff -1}
\qff =\qff
I_t\dff U_t^*
\qff =\qff
P_t^*\dff U_t^*
\qff =\qff
(U_t\dff P_t)^*$\nsp.\qss
Therefore the norm continuity of $I_t\dff U_t^{\dff -1}$ 
follows from the norm continuity of $U_t\dff P_t$\nsp.\qss
The proof for $I'_t\dff V_t^{\dff -1}$
is completely similar.\qss

\myuppar{Construction of\dss Dixmier--Douady--like families.}
Let us fix an orthonormal\sss basis $e_1\fff,\dff e_2\fff,\dff \ldots$ of\dss $H$\nnsp.\qss
Let $H_{\dff 1/2}$ be a closed subspace of\sss $H$ of\sss infinite dimension and codimension
and $H'_{\fff 1/2}$ be its orthogonal complement.\qss
Let $u\colon H\ttoo H$ and $v\colon H\ttoo H$
be some isometric embeddings with the images $H_{\dff 1/2}$ and $H'_{\fff 1/2}$
respectively.\qss
Arguing as in the proof\sss of\dss Theorem\qss \ref{t13},\qss
but using the path connectedness of\sss the unit sphere instead of\dss
Kakutani's\dss theorem,\qss
we can construct a norm continuous family $w_t$\nsp, $t\in [1/2,\dff 1)$ of\sss unitary operators in $H$
such that $w_{\dff 1/2}\qff =\qff \id_{\dff H}$ and $w_t(u(e_i))\qff =\qff e_i$ for
$t\geq 1-1/(i+2)$\nnsp.\qss
Then $u_t\qff =\qff w_t\dff u$\nnsp, $t\in [1/2,\dff 1)$
is a norm continuous family of\sss isometric embeddings $H\ttoo H$
such that $u_{\dff 1/2}\qff =\qff u$ and $u_t(e_i)\qff =\qff e_i$ for
$t\geq 1-1/(i+2)$\nnsp.\qss

Lemma\qss \ref{l10}\qss (more precisely,\qss its first\sss part)\qss
implies that $u_t$ strongly converges to $\id_{\dff H}$ when $t\ttoo 1$\nnsp.\qss
Let us set $u_1\qff =\qff \id_{\dff H}$\nsp.\qss
The norm continuity of\sss the family $u_t$\nsp, $t\in [1/2,\dff 1)$
implies that there exist a norm continuous family $v_t$\nsp, $t\in [1/2,\dff 1)$
such that the image of $v_t$ is equal\sss to the orthogonal complement 
of\sss the image of $u_t$ for every $ t\in [1/2,\dff 1)$\nnsp.\qss
Alternatively,\qss one can simply set\sss 
$v_t\qff =\qff w_t\dff v$ for every $t\in [1/2,\dff 1)$\nnsp.\qss

Similarly,\qss one can construct a norm continuous family 
$v_t$\nsp, $t\in (0,\dff 1/2]$ of\sss isometric embeddings
$H\ttoo H$ such that $v_{\dff 1/2}\qff =\qff v$
and $v_t$ strongly converges to $\id_{\dff H}$ when $t\ttoo 0$\nnsp,\qss
and\sss then set $v_0\qff =\qff \id_{\dff H}$ and choose a norm continuous family
$u_t$\nsp, $t\in (0,\dff 1/2]$ such that the image of $u_t$ is equal\sss to the orthogonal complement 
of\sss the image of $v_t$ for every $ t\in (0,\dff 1/2]$\nnsp.\qss

Now,\qss let us set  $H_{\fff 0}=\dff 0$\nnsp, $H_{\fff 1}=\dff H$ 
and define the subspaces $H_t$ and $H'_t$ for $t\in (0,\dff 1)$ of\sss $H$
as the images of\sss $u_t$ and $v_t$ respectively.\qss
Then $H'_t\qff =\qff H\ominus H_t$ is the orthogonal complement of $H_t$\nsp.\qss
Let\sss us define $P_t\dff,\qff I_t\dff,\qff Q_t\dff,\qff I'_t$ as above.\qss
The maps $u_t$ and $v_t$ induce isometric isomorphisms
$H\ttoo H_t$ and $H\ttoo H'_t$ respectively.\qss
Let $U_t$ and $V_t$ be,\qss respectively,\qss their inverses.\qss
Then $u_t\qff =\qff I_t\dff U_t^{\dff -1}$ and
$v_t\qff =\qff I'_t\dff V_t^{\dff -1}$\dnsp.\qss
The continuity properties of\sss these operators hold\sss by\sss the construction.\qss
When $t\ttoo 1$\nnsp,\qss the operator $I_t\dff P_t$ strongly converges to $\id_{\dff H}$\nsp.\qss
Indeed,\qss for every $x\in H$ we have
$\norm{I_t\dff P_t(x)\qff -\qff x}\qff \leq\qff \norm{u_t(x)\qff -\qff x}$
because $P_t(x)$ is\sss the closest to $x$ point of\sss the image of $u_t$\nsp.\qss
Since $u_t(x)$ converges to $x$\nnsp,\qss
this implies that $I_t\dff P_t(x)$ also converges to $x$\nnsp.\qss
Hence $I_t\dff P_t$ strongly converges to $\id_{\dff H}$\nsp.\qss
Since $I'_t\dff Q_t\qff =\qff \id_{\dff H}\dff -\qff I_t\dff P_t$\nsp,\qss
it\sss follows that $I'_t\dff Q_t$ strongly converges to $0$ when $t\ttoo 1$\nnsp.\qss
Similarly,\qss when $t\ttoo 0$\nnsp,\qss the operators $I'_t\dff Q_t$ and $I_t\dff P_t$ 
strongly converge to $\id_{\dff H}$ and $0$ respectively.\qss
We see that our families have all\sss promised properties.\qss

\myuppar{Another proof\dss of\dss the contractibility of\sss $\mathcal{U}(H)$\nnsp.}
For  $U\in\dff \mathcal{U}(H)$ and $t\in (0,\dff 1]$ let
\[
\quad
\Phi(U,\dff t)
\qff =\qff
(\id_{\dff H}\dff -\dff I_t\dff P_t)
\qff +\qff
I_t\dff U_t^{\dff -1}\dff U\dff U_t\dff P_t
\qff.
\]
Let $\Phi(U,\dff 0)=\id_{\dff H}$\nsp.\qss
Then $\Phi(U,\dff 1)=U$
and the operators 
$\Phi(U,\dff t)$
strongly continuously depend on $(U,\dff t)$ for $t\in (0,\dff 1]$\nnsp.\qss
This follows from the fact that the operators
$I_t\dff P_t$\nsp,\dss $U_t\dff P_t$\nsp,\qss and $I_t\dff U_t^{\dff -1}$
norm continuously depend on $t$ for $t\in (0,\dff 1)$ and
strongly continuously depend on $t$ at $t=1$\nnsp.\qss
At the same time when $t\ttoo 0$
the operator $\id_{\dff H}\dff -\dff I_t\dff P_t$ strongly converges to $\id_{\dff H}$ and
$\norm{\fff I_t\dff U_t^{\dff -1}\dff U\dff U_t\dff P_t(f)}
\qff =\qff
\norm{\fff P_t(f)}\ttoo 0$
for every $f$ because $I_t\dff P_t$ strongly converges to $0$\nnsp.\qss
It follows that $\Phi(U,\dff t)$ strongly continuously depends on $t$
at $t=0$\nnsp.\qss
Hence $\Phi(U,\dff t)$ defines a homotopy
between the identity map of $\mathcal{U}(H)$ and the map taking
the whole $\mathcal{U}(H)$ to $\id_{\dff H}$\nsp.\qss

\myuppar{Remarks.}
The above proof follows the proof of\dss Dixmier--Douady\qss \cite{dd}.\qss
See\qss \cite{dd},\qss Lemma\qss 3.\qss
Of course,\qss Dixmier--Douady\dss use different families $P_t\dff,\dff U_t$\nsp.\qss
Their families are only strongly continuous for $t\in (0,\dff 1)$\nnsp,\qss
in contrast with the above ones,\qss which are norm continuous there.\qss

The readers who are going to compare the above proof with\qss \cite{dd}\qss
should keep in mind that\dss Dixmier--Douady do not distinguish
between $P_t$ and $I_t\dff P_t$\nsp,\qss
as also between $U_t^{\dff -1}$ and $I_t\dff U_t^{\dff -1}$\nnsp.\qss

\myuppar{Compact operators.}
The arguments in\sss the second proof of\sss the contractibility of $\mathcal{U}(H)$
can be modified\sss to apply to compact operators.\qss
Let $K(H)$ be the space of compact operators in $H$ with the norm topology.\qss
Let us consider the map $K(H)\times (0,\dff 1]\ttoo K(H)$ defined by
\[
\quad
(A,\dff t)
\longmapsto
I_t\dff U_t^{\dff -1}\dff A\dff U_t\dff P_t
\qff.
\]
We will\sss prove that it is continuous in the norm topology, including at $t=1$\nnsp.\qss
With our choice of operators\sss $P_t\dff,\qff U_t$
the norm continuity on $K(H)\times (0,\dff 1)$ is trivial.\qss
Therefore it\sss is sufficient to prove the continuity at each point of the form $(A,\dff 1)$\nnsp.\qss
In order to do this,\qss we need to estimate
$\norm{I_t\dff U_t^{\dff -1}\dff B\dff U_t\dff P_t
\qff -\qff
A}$\nnsp,\qss
where $B\in K(H)$\nnsp,\qss in terms of $\norm{B-A}$ and $t$\nnsp.\qss
Clearly,
\[
\quad
\norm{I_t\dff U_t^{-1}\dff B\dff U_t\dff P_t\qff -\qff I_t\dff U_t^{-1}\dff A\dff U_t\dff P_t}
\qff \leq\qff
\norm{B-A}
\qff.
\]
It\sss follows that\sss it\sss is sufficient to estimate
$\norm{I_t\dff U_t^{\dff -1}\dff A\dff U_t\dff P_t
\qff -\qff
A}$\nnsp.\qss
Recalling that
$I_t\dff U_t^{\dff -1}
\qff =\qff
(U_t\dff P_t)^*$
and using the fact that
$\norm{D^*}=\norm{D}$
for every bounded operator $D$\dnsp,\qss
we see that
\[
\quad
\norm{I_t\dff U_t^{-1}\dff A\dff U_t\dff P_t\qff -\qff A}
\qff \leq\qff
\norm{I_t\dff U_t^{-1}\dff A\dff U_t\dff P_t\qff -\qff A\dff U_t\dff P_t}
\qff +\qff
\norm{A\dff U_t\dff P_t\qff -\qff A}
\]

\vspace{-36pt}
\[
\quad
=\qff
\norm{(\fff I_t\dff U_t^{-1}\qff -\qff \id_{\dff H})\dff A\dff U_t\dff P_t}
\qff +\qff
\norm{((U_t\dff P_t)^*\qff -\qff \id_{\dff H})\dff A^*}
\]

\vspace{-36pt}
\[
\quad
=\qff
\norm{(\fff I_t\dff U_t^{-1}\qff -\qff \id_{\dff H})\dff A\dff U_t\dff P_t}
\qff +\qff
\norm{(\fff I_t\dff U_t^{\dff -1}\qff -\qff \id_{\dff H})\dff A^*}
\]
Since $A$ is compact, the images of\sss the unit ball under the operators
$A\dff U_t\dff P_t$ and $A^*$ are contained in a compact set independent of $t$\nnsp.\qss
Since the operators $I_t\dff U_t^{\dff -1}\qff -\qff \id_{\dff H}$ tend to $0$ when $t\ttoo 1$
in the strong topology and have the norm $\leq 2$\nnsp,\qss 
they also tend to $0$ in the compact-open topology.\qss
Therefore both summands in the last displayed formula tend to $0$ when $t\ttoo 1$\nnsp.\qss
This provides the needed estimate and\sss hence proves the continuity.\qss

Similarly,\qss the map $K(H)\times [0,\dff 1)\ttoo K(H)$
defined by
$(A,\dff t)
\longmapsto
I'_t\dff V_t^{\dff -1}\dff A\dff V_t\dff Q_t$
is also norm continuous.\qss
This leads to a norm continuous map
$K(H)\times [0,\dff 1]\times K(H) \ttoo K(H)$
defined by
\[
\quad
(A,\dff t\dff, B)
\longmapsto
t\left(\fff I_t\dff U_t^{\dff -1}\dff A\dff U_t\dff P_t\right)
\qff +\qff
(1-t)  \left(\fff I'_t\dff V_t^{\dff -1}\dff B\dff V_t\dff Q_t\right)
\qff.
\]
The norm continuity of the fist summand at $t=0$ is ensured by the factor $t$\nnsp,\qss
and of the second summand at $t=1$ by the factor $1-t$\nnsp.\qss
One can think about this map as a canonical\sss homotopy between
arbitrary two compact operators $A,\dff B$ continuously depending on $A,\dff B$\nnsp.\qss

\myuppar{Remark.}
The above proof\sss of\sss the continuity is an adaptation of\sss the proof of\dss
Proposition\qss A.1.1\qss in\dss Atiyah--Segal\qss \cite{as}.\qss
The map
$K(H)\times [0,\dff 1]\times K(H) \ttoo K(H)$
is a minor generalization of a construction from\qss \cite{p}.\qss
See\qss \cite{p},\qss the proof of\dss Proposition\dss 4.1.\qss

\begin{flushright}

March\qss 21,\oss 2025

\end{flushright}

N.I.\dff:\quad https\halfff:/\!/nikolaivivanov.com;\quad
E-mail\halfff:\oss nikolai.v.ivanov{\fff}@{\dff}icloud.com,\oss ivanov{\fff}@{\dff}msu.edu;\\
Department\sss of\qss Mathematics,\oss Michigan\sss State\sss University

M.P.\dff:\quad 
https\halfff:/\!/marina-p.info;\quad
E-mail\halfff:\oss marina.p@technion.ac.il;\quad\\
Department\sss of\qss Mathematics,\qss University of\dss Haifa\qss (Israel);\qss\\
Department\sss of\qss Mathematics,\qss Technion -- Israel\dss Institute of\dss Technology

\end{document}